\documentclass[12pt,a4paper]{article}
\usepackage{amssymb}

\usepackage{graphicx,amssymb,amsfonts,epsfig,amsthm,a4,amsmath,url}
\usepackage[latin1]{inputenc}

\newtheorem{Thmm}{Theorem}

\newtheorem{Propp}[Thmm]{Proposition}
\newtheorem{Corr}[Thmm]{Corollary}

\newtheorem{Thm}{Theorem}[section]
\newtheorem{Prop}[Thm]{Proposition}
\newtheorem{Lem}[Thm]{Lemma}
\newtheorem{Cor}[Thm]{Corollary}
\theoremstyle{definition}
\newtheorem{Def}[Thm]{Definition}

\newtheorem*{Ack}{Acknowledgement}
\newtheorem{Rem}[Thm]{Remark}
\newtheorem{Ex}[Thm]{Example}

\newtheorem{quest}{Question}

\newcommand{\Z}{\mathbf{Z}}

\newcommand{\N}{\mathbf{N}}
\newcommand{\R}{\mathbf{R}}

\newcommand{\HH}{\mathcal{H}}

\newcommand{\eps}{\varepsilon}

\title{Isometric group actions on Hilbert spaces:
structure of orbits}
\author{Yves de Cornulier, Romain Tessera, Alain Valette}
\date{\today}

\begin{document}

\baselineskip=16pt

\maketitle

\begin{abstract}
Our main result is that a finitely generated nilpotent group has
no isometric action on an infinite-dimensional Hilbert space with
dense orbits. In contrast, we construct such an action with a
finitely generated metabelian group.

\medskip

    \hfill\break
\noindent {\sl Mathematics Subject Classification:} Primary 22D10;
Secondary 43A35, 20F69. \hfill\break {\sl Key words and Phrases:}
Affine actions, Hilbert spaces, minimal actions, nilpotent groups.
\end{abstract}

\section{Introduction}

The study of isometric actions of groups on affine Hilbert spaces
has, in recent years, found applications ranging from the
$K$-theory of $C^*$-algebras \cite{HiKa}, to rigidity theory
\cite{ShaInv} and geometric group theory \cite{Shalom04,CoTeVa1}.
This renewed interest motivates the following general problem:
{\it How can a given group act by isometries on an affine Hilbert
space?}

This paper is a sequel to \cite{CoTeVa1}, but can be read
independently. In \cite{CoTeVa1}, given an isometric action of
a finitely generated group $G$ on a Hilbert space
$\alpha:G\to\textnormal{Isom}(\mathcal{H})$, we focused on the
growth of the function $g\mapsto\alpha(g)(0)$. Here the emphasis
is on the structure of orbits.

We will mainly focus on actions of nilpotent groups. Let us begin
by a simple example: every isometric action of $\Z$ on a Euclidean
space is the direct sum of an action with a fixed point and an
action by translations. This actually remains true for general
nilpotent groups. The situation becomes more
subtle when we study actions on infinite-dimensional Hilbert
spaces. However, something remains from the finite-dimensional
case.

We say that a convex subset of a Hilbert space is \textit{locally
bounded} if its intersection with any finite-dimensional subspace
is bounded. The main result of the paper is the following theorem, proved 
in \S\ref{Sec:nonex}.

\begin{Thmm}Let $G$ be a nilpotent topological group. Let $G$ act 
isometrically on a
Hilbert space $\mathcal{H}$, with linear part $\pi$. Let
$\mathcal{O}$ be an orbit under this action. Then there
exist\begin{itemize}\item a closed subspace $T$ of $\HH$ (the
``translation part"), contained in the subspace of invariant vectors of 
$\pi$,
and\item a closed, locally bounded convex subset $U$ of the
orthogonal subspace $T^\bot$,\end{itemize} such that $\mathcal{O}$
is contained in $T\times U$.\label{main_orbits}
\end{Thmm}

We owe the following general question to A.~Navas: which locally
compact groups have an isometric action on an infinite-dimensional
separable Hilbert space with dense orbits (i.e. a minimal action)?
Theorem \ref{main_orbits} allows us to provide a negative answer
in the case of finitely generated nilpotent groups.

\begin{Corr}\textnormal{(see Corollary \ref{novirtuallynilpotent})}
A compactly generated, nilpotent-by-compact locally compact group
does not admit any affine isometric action with dense orbits on an
infinite-dimensional Hilbert space.\label{cor_dense_nilp}
\end{Corr}
%

In the course of our proof, we introduce the following new
definition: a unitary or orthogonal representation $\pi$ of a
group is \textit{strongly cohomological} if $H^1(G,\rho)\neq 0$ for
every nonzero subrepresentation $\rho\le\pi$. It is easy to observe that 
the linear part of
an affine isometric action with dense orbits is strongly
cohomological. The main non-trivial step in the proof of Theorem 
\ref{main_orbits} is the following result.

\begin{Propp}\textnormal{(see Proposition \ref{cor:nilp_str_coh_triv})}
Let $\pi$ be an orthogonal or unitary representation of a second
countable, nilpotent locally compact group $G$. Suppose that $\pi$
is strongly cohomological. Then $\pi$ is a trivial representation.
\end{Propp}

Another case for which we answer negatively Navas' question is the 
following.

\begin{Thmm}\textnormal{(see Theorem \ref{thm:nodense_semisimple})}
Let $G$ be a connected semisimple Lie group. Then $G$ has no
isometric action on a nonzero Hilbert space with dense orbits.
\end{Thmm}

It is not clear how Theorem \ref{main_orbits} and Corollary 
\ref{cor_dense_nilp} can be generalized, in view
of the following example.

\begin{Propp}\textnormal{(see Proposition
\ref{Prop:exist_metab_dense_Hilbert})} There exists a finitely
generated metabelian group admitting an affine isometric action
with dense orbits on an infinite-dimensional separable 
Hilbert space.
\end{Propp}

Another construction provides

\begin{Propp}\textnormal{(see Proposition \ref{minaction-almostfixed})}
There exists a countable group admitting an affine isometric
action with dense orbits on an infinite-dimensional Hilbert space,
in such a way that every finitely generated subgroup has a fixed
point.
\end{Propp}

\begin{Ack}We thank A.~Navas for useful discussions and encouragement.\end{Ack}

\section{Existence results}

Here is a first positive result regarding Navas' question.

\begin{Prop}
There exists an isometric action of a metabelian 3-generator group
on $\ell^2_{\R}(\Z)$, all 
of whose
orbits are dense.\label{Prop:exist_metab_dense_Hilbert}
\end{Prop}

\begin{proof} Observe that $\Z[\sqrt{2}]$ acts on $\R$ by translations, 
with
dense orbits. So the free abelian group of countable rank
$\Z[\sqrt{2}]^{(\Z)}$ acts by translations, with dense orbits, on
$\ell^2_\R(\Z)$. Observe now that the latter action extends to the
wreath product
$\Z[\sqrt{2}]\wr\Z=\Z[\sqrt{2}]^{(\Z)}\rtimes\penalty10000\Z$,
where $\Z$ acts on $\ell^2_\R(\Z)$ by the shift. That wreath
product is metabelian, with 3 generators.\end{proof}

\begin{Cor}
There exists an isometric action of a free group of finite rank on
a Hilbert space, with dense orbits.\qed
\end{Cor}

Recall that an isometric action
$\alpha:G\to\textnormal{Isom}(\mathcal{H})$ \textit{almost has
fixed points} if for every $\eps>0$ and every compact subset
$K\subset G$ there exists $v\in\mathcal{H}$ such that $\sup_{g\in
K}\|v-\alpha(g)v\|\le\eps$.

In the example given by Proposition
\ref{Prop:exist_metab_dense_Hilbert}, the given isometric action
clearly does not almost have fixed points, i.e. it defines a
nonzero element in reduced 1-cohomology. The next result shows
that this is not always the case.

\begin{Prop}\label{minaction-almostfixed} There exists a countable
group $\Gamma$ with an affine isometric action $\alpha$ on an 
infinite-dimensional
Hilbert space, such that $\alpha$ has dense orbits, and every
finitely generated subgroup of $\Gamma$ has a fixed point. In
particular, the action almost has fixed points.
\end{Prop}

\begin{proof}We first construct an uncountable group $G$ and an affine
isometric action of $G$ having dense orbits and almost having fixed
points.

In $\mathcal{H}\,=\,\ell^{2}_{\R}(\N)$, let $A_n$ be the affine subspace 
defined by the equations
$$x_0 =1,\,x_1 =1,...,\,x_n =1,$$
and let $G_n$ be the pointwise stabilizer of $A_n$ in the isometry
group of $\mathcal{H}$. Let $G$ be the union of the $G_n$'s. View
$G$ as a discrete group.

It is clear that $G$ almost has fixed points in $\mathcal{H}$,
since any finite subset of $G$ has a fixed point. Let us prove
that $G$ has dense orbits.

\noindent {\bf Claim 1.} For all $x,y\in\mathcal{H}$, we have
$\lim_{n\rightarrow\infty} |d(x,A_n)-d(y,A_n)|=0$.

By density, it is enough to prove Claim 1 when $x,y$ are finitely
supported in $\ell^{2}_{\R}(\N)$. Take
$x=(x_0,x_1,...,x_k,0,0,...)$ and choose $n>k$. Then
$$d(x,A_n)^2 \,=\,\sum_{j=0}^{k}(x_j -1)^2 + \sum_{j=k+1}^{n}1^2
\,=\,n+1 -2\sum_{j=0}^{k}x_j + \sum_{j=0}^{k}x_j^2,$$ so that
$d(x,A_n)=\sqrt{n} + O(\frac{1}{\sqrt{n}})$, which proves Claim 1.

\medskip

Denote by $p_n$ the projection onto the closed convex set $A_n$,
namely $p_n(x_0,x_1,\dots)=(1,1,\dots,1,x_{n+1},x_{n+2},\dots)$.

\noindent {\bf Claim 2.} For all $x,y\in\mathcal{H}$, we have
$\lim_{n\rightarrow\infty} \|p_n(x)-p_n(y)\|=0$.

This is a straightforward computation.

\medskip

\noindent {\bf Claim 3.} $G$ has dense orbits in $\mathcal{H}$.

Observe that two points $x,y\in\mathcal{H}$ are in the same
$G_n$-orbit if and only if $d(x,A_n)=d(y,A_n)$ and
$p_n(x)=p_n(y)$. Fix $x_0, z\in\mathcal{H}$. We want to show that
$\lim_{n\rightarrow\infty} d(G_n x_0,z)=0$. So fix $\eps
>0$. By the second claim, for some $n_0$, $\|p_n(x_0)-p_n(z)\|\le\eps/2$ whenever
$n\ge n_0$. Set
$$W=\{x\in\mathcal{H}:\,p_n(x)=p_n(z)\};$$ this is the
orthogonal affine subspace of $A_n$ passing through $z$. Then
$y_0=x_0+(p_n(z)-p_n(x_0))\in W$. By the first claim, there exists
$n_1\ge n_0$ such that $|d(y_0,A_n)-d(z,A_n)|\le\eps/2$ for every
$n\ge n_1$. Therefore there exists $y\in W$ such that $\|y-z\|\le
\eps/2$ and $d(y,A_n)=d(y_0,A_n)=d(x_0,A_n)$. By the previous
observation, there exists $g\in G_n$ such that $y=gy_0$. Then
$$d(gx_0,z)\le d(gx_0,gy_0)+d(gy_0,z)\le\eps,$$
so that $d(G_nx_0,z)\le\eps$ for every $n\ge n_1$, proving the
last claim.

\medskip

Using separability of $\mathcal{H}$, it is now easy to construct a
countable subgroup $\Gamma$ of $G$ also having dense orbits on
$\mathcal{H}$.\end{proof}

\begin{quest} Does there exist an affine isometric action
of a {\it finitely generated} group on a Hilbert space, having
dense orbits and almost having fixed points?\end{quest}

\section{Cohomology of unitary representations of nilpotent
groups}

Our non-existence results concerning nilpotent
groups will be based on the following study of their unitary
representations.

\begin{Def}
If $G$ is a topological group and $\pi$ a unitary representation,
we say that $\pi$ is {\em strongly cohomological} if every nonzero
subrepresentation of $\pi$ has nonzero first cohomology.
\end{Def}

The following lemma is Proposition 3.1 in Chapitre III of
\cite{Guich}.

\begin{Lem}
Let $\pi$ be a unitary representation of a topological group $G$. Let $z$ 
be a central element 
of $G$. Suppose that $1-\pi(z)$ has a bounded inverse (equivalently, 1 
does not belong to the spectrum of $\pi(z)$). Then 
$H^1(G,\pi)=0$.\label{Lem:Gui31}
\end{Lem}
\begin{proof} Let $b\in Z^1(G,\pi)$ be a 1-cocycle; we prove that $b$ is bounded. 
If $g\in G$, expanding the equality $b(gz)=b(zg)$, we
obtain that $(1-\pi(z))b(g)$ is bounded by $2\|b(z)\|$, so that
$b$ is bounded by $2\|(1-\nolinebreak\pi(z))^{-1}\|\|b(z)\|$.
\end{proof}

\begin{Lem}\label{trivialonZ(G)}
Let $G$ be a locally compact, second countable group, and $\pi$ a
strongly cohomological unitary representation. Then $\pi$ is trivial on
the centre $Z(G)$.
\end{Lem}

\begin{proof}Fix $z\in Z(G)$. As $G$ is second countable, we may write
$\pi=\int_{\hat{G}}^{\oplus}\rho\,d\mu(\rho)$, a disintegration of
$\pi$ as a direct integral of irreducible representations. Let
$\chi:\hat{G}\rightarrow S^1:\rho\mapsto\rho(z)$ be the continuous
map given by the value of the central character of $\rho$ on $z$.
For $\eps>0$, set
$X_{\eps}=\{\rho\in\hat{G}:|\chi(\rho)-1|>\eps\}$ and
$\pi_{\eps}=\int_{X_{\eps}}^{\oplus}\rho\,d\mu(\rho)$, so that
$\pi_{\eps}$ is a subrepresentation of $\pi$. Since
$|\rho(z)-1|^{-1}<\eps^{-1}$ for $\rho\in X_{\eps}$, the operator
$$(\pi_{\eps}(z)-1)^{-1}\,=\,\int_{X_{\eps}}^{\oplus}(\rho(z)-1)^{-1}
\,d\mu(\rho)$$ is bounded. We can now apply Lemma
\ref{Lem:Gui31}
to conclude that $H^1(G,\pi_{\eps})=0$. By definition, this means
that $\pi_{\eps}$ is the zero subrepresentation, meaning that the
spectral measure $\mu$ is supported in $\hat{G}-X_{\eps}$. As this holds
for every $\eps>0$, we see that $\mu$ is supported in
$\{\rho\in\hat{G}:\,\rho(z)=1\}$, to the effect that
$\pi(z)=1$.\end{proof}

\begin{Prop}
Let $G$ be a topological group, and $\pi$ a unitary representation
of $G$. Suppose that $\overline{H^1}(G,\pi)\neq 0$. Then $\pi$ has
a nonzero subrepresentation that is strongly
cohomological.\label{Prop:H1_bar_vs_str_coh}
\end{Prop}

\begin{proof}Suppose the contrary. Then, by a standard application of
Zorn's Lemma, $\pi$ decomposes as a direct sum
$\pi=\bigoplus_{i\in I}^{}\pi_i$, where $H^1(G,\pi_i)=0$ for every
$i\in I$, so that $\overline{H^1}(G,\pi)=0$ by Proposition 2.6 in
Chapitre III of \cite{Guich}.\end{proof}

\begin{Rem}
The converse is false, even for finitely generated groups: indeed,
it is easy to check (see \cite{Gui72}) that every nonzero unitary
representation of the free group $F_2$ has non-vanishing $H^1$, so
that every unitary representation of $F_2$ is strongly
cohomological. But it turns out that $F_2$ has an irreducible
representation $\pi$ such that $\overline{H^1}(F_2,\pi)=0$ (see
Proposition 2.4 in \cite{MarVal}).
\end{Rem}

\begin{Cor}\label{moduloZ(G)}
Let $G$ be a locally compact, second countable group, and let
$\pi$ be a unitary representation of $G$ without invariant
vectors. Write $\pi=\pi_0\oplus\pi_1$, where $\pi_1$ consists of
the $Z(G)$-invariant vectors. Then
\begin{enumerate}
    \item  [(1)] $\pi_0$ does not contain any
nonzero strongly cohomological subrepresentation; in particular,
$\overline{H^1}(G,\pi_0)=0$;

    \item  [(2)] every 1-cocycle of $\pi_1$
vanishes on $Z(G)$, so that \allowbreak$H^1(G,\pi_1)\simeq
H^1(G/Z(G),\pi_1)$.
\end{enumerate}
\end{Cor}

\begin{proof}(1) follows by combining Lemma \ref{trivialonZ(G)} and
Proposition \ref{Prop:H1_bar_vs_str_coh}. For (2), we use the idea
of proof of Theorem~3.1 in \cite{ShaInv}: if $b\in
Z^1(G,\pi_{1})$, then for every $g\in G,\,z\in Z(G)$,
$$\pi_{1}(g)b(z)+b(g)\,=\,b(gz)\,=\,b(zg)\,=\,b(g)+b(z)$$
as $\pi_{1}(z)=1$. So $\pi_{1}(g)b(z)=b(z)$; this forces $b(z)=0$
as $\pi$ has no $G$-invariant vector. So $b$ factors through
$G/Z(G)$.\end{proof}

Observe that Corollary \ref{moduloZ(G)} provides a new proof of
Shalom's Corollary 3.7 in \cite{ShaInv}: under the same
assumptions, every cocycle in $Z^1(G,\pi)$ is almost cohomologous
to a cocycle factoring through $G/Z(G)$ and taking values in a
subrepresentation factoring through $G/Z(G)$. From Corollary 
\ref{moduloZ(G)} we also immediately deduce

\begin{Cor}\label{cohnilp}
Let $G$ be a locally compact, second countable, nilpotent group,
and let $\pi$ be a unitary representation of $G$ without invariant
vectors. Let $(Z_i)$ be the ascending central series of $G$
($Z_0=\{1\}$, and $Z_i$ is the centre modulo $Z_{i-1}$). Let
$\sigma_i$ denote the subrepresentation of $G$ on the space of $Z_i$-invariant
vectors, and finally let $\pi_i$ be
the orthogonal of $\sigma_{i+1}$ in $\sigma_{i}$, so that
$\pi=\bigoplus\pi_i$. 

Then $H^1(G,\pi_i)\simeq H^1(G/Z_i,\pi_i)$ for all $i$, and $\pi$
has no nonzero strongly cohomological subrepresentation. In particular,
$\overline{H^1}(G,\pi)=0$.
\hfill$\square$
\end{Cor}

Note that the latter statement, namely $\overline{H^1}(G,\pi)=0$, is a 
result of Guichardet
\cite[Th\'e\-or\`eme~7]{Gui72}, which can be stated as: $G$ has
Property $H_{T}$ (i.e. every unitary representation with
non-vanishing reduced 1-cohomology contains the trivial
representation).

\begin{Def}
We say that a locally compact group $G$ has Property $H_{CT}$ if
every strongly cohomological unitary representation of $G$ is
trivial.\end{Def}

It is a straightforward verification that this is equivalent to:
every strongly cohomological \textit{orthogonal} representation of
$G$ is trivial. This will be useful in the next paragraph since we
will deal with orthogonal rather than unitary representations.
The following proposition is contained in Corollary \ref{cohnilp}.

\begin{Prop}
If $G$ is a locally compact, second countable nilpotent group,
then $G$ has Property $H_{CT}$.\qed\label{cor:nilp_str_coh_triv}
\end{Prop}

As a corollary of Proposition \ref{Prop:H1_bar_vs_str_coh}, Property 
$H_{CT}$ implies Property $H_T$. However the converse is not true, as 
shown by the following example.

\begin{Ex} Let $G$ be the full affine group of the real line. The dual 
$\hat{G}$ (i.e. the space of unitary irreducible representations of $G$ 
with the Fell-Jacobson topology) was described in \cite{Fell}: it consists 
of two copies of the real line (corresponding to one-dimensional 
representations, i.e. characters) plus one point $\{\sigma\}$ which is 
both open and dense. The only irreducible representation with 
non-vanishing reduced 1-cohomology is the trivial representation $1_G$, so 
that $G$ has Property $H_T$; on the other hand, since $\sigma$ weakly 
contains $1_G$, one has $H^1(G,\sigma)\neq 0$ by 
\cite[Th\'eor\`eme~1]{Gui72}. So $\sigma$ is strongly cohomological, meaning 
that $G$ 
does not have Property $H_{CT}$.\end{Ex}

\section{Non-existence results}\label{Sec:nonex}

\begin{Def} 1) We say that a subset $Y$ of a metric space
$(X,d)$ is {\it coarsely dense} if there exists $C\geq 0$ such
that, for every $x\in X$,
$$d(x,Y)\leq C.$$

2) We say that a subset $Y$ of a Hilbert space $\mathcal{H}$ is
\textit{enveloping} if its closed convex hull is all of
$\mathcal{H}$.
\end{Def}

Observe that every dense subset of a metric space is coarsely
dense. Besides, in a Hilbert space $\mathcal{H}$, every coarsely
dense subset $Y$ is enveloping. Indeed, suppose that $Y$ is
contained in a closed, convex proper subset $X$ of $\mathcal{H}$.
Consider $v\notin X$ and let $y$ denote its projection on $X$
(excluding the trivial case $Y=\emptyset$). Then, for every
$\lambda\ge 0$, we have $d(y+\lambda(v-y),Y)\ge
d(y+\lambda(v-y),X)=\lambda$, which is unbounded, so that $Y$ is
not coarsely dense.

\begin{Ex}
In $\ell^2_\R(\Z)$, let $X$ denote the subset of elements with integer
coefficients. Then $X$ is enveloping: indeed, its intersection
with the subspace $V_n=\ell^2_\R(\{-n,\dots,n\})$ is coarsely
dense, hence enveloping in $V_n$, and the increasing union
$\bigcup V_n$ is dense in $\ell^2_\R(\Z)$. But $X$ is not coarsely
dense: indeed, for every $n\ge 0$, the element
$\frac12\mathbf{1}_{\{1,\dots,4n\}}$ is at distance $\sqrt{n}$ to
$X$.

Note that $X$ is the orbit of $0$ for the natural action of the
wreath product $\Z\wr \Z=\Z^{(\Z)}\rtimes\Z$ on $\ell^2_\R(\Z)$,
where $\Z^{(\Z)}$ acts by translations and the factor $\Z$ acts by
shifting (compare to the example in the proof of Proposition
\ref{Prop:exist_metab_dense_Hilbert}).
\end{Ex}

\begin{Lem}
Let $G$ be a topological group and $\pi$ an orthogonal
representation, admitting a 1-cocycle $b$ with enveloping orbits.
Then $\pi$ is strongly cohomological.\label{Lem:dense_strong_coh}
\end{Lem}

\begin{proof}If $\sigma$ is a nonzero subrepresentation of $\pi$, let
$b_{\sigma}$ be the orthogonal projection of $b$ on
$\HH_{\sigma}$, so that $b_{\sigma}\in Z^1(G,\sigma)$. Then
$b_{\sigma}(G)$ is enveloping in $\HH_{\sigma}$, in particular
$b_{\sigma}$ is unbounded. So $b_{\sigma}$ defines a nonzero
class in $H^1(G,\sigma)$.\end{proof}

\begin{Thm}
Let $G$ be a topological group with Property $H_{CT}$. Let $G$
act isometrically on a Hilbert space $\mathcal{H}$, with linear
part $\pi$. Let $\mathcal{O}$ be an orbit under this action. Then
there exist\begin{itemize}\item a subspace $T$ of $\HH$, contained
in $\mathcal{H}^{\pi(G)}$, and\item a closed, locally bounded
convex subset $U$ of $T^\bot$,\end{itemize} such that
$\mathcal{O}$ is contained in $T\times U$.\label{thm:orbits_hct}
\end{Thm}
\begin{proof}
We immediately reduce to the case when $\pi$ has no invariant
vectors, so that we must prove that the closed convex hull $U$ of
$\mathcal{O}$ is locally bounded.

Observe that a convex subset of a Hilbert space is locally bounded
if and only if it contains no affine half-line. Thus denote by
$\mathcal{D}$ the set of affine half-lines contained in $U$, and
suppose by contradiction that $\mathcal{D}\neq\emptyset$. Denote
by $\mathcal{D}_0$ the corresponding set of linear half-lines
(where the linear half-line corresponding to a half-line $x+\R_+v$
is simply $\R_+v$). Then $\mathcal{D}_0$ is invariant under the
linear action $\pi$ of $G$. Let $W$ be the closed subspace of
$\mathcal{H}$ generated by all the half-lines in $\mathcal{D}_0$,
and denote by $\sigma$ the corresponding subrepresentation. By
assumption, $\sigma$ is nonzero.

We claim that $\sigma$ is strongly cohomological, contradicting
that $\pi$ has no invariant vectors along with the $H_{CT}$
assumption. Let $\rho$ be a nonzero subrepresentation of
$\sigma$. Then by the definition of $W$, there exists an half-line
of $U$ which projects injectively into the subspace of $\rho$.
Thus $H^1(G,\rho)\neq 0$, proving the claim, and ending the
proof.\end{proof}

\begin{proof}[Proof of Theorem \ref{main_orbits}]
We can suppose that $\pi$ has no invariant vectors. Suppose that the 
convex hull of $\alpha(G)(0)$ is not locally bounded. Then it contains a 
half-line $D=x+\R_+v$. Let $(x_n)$ be an unbounded sequence in $D$. Every 
$x_n$ is a convex combination of elements of the form $\alpha(g)(0)$, 
where $g$ ranges over a finite subset $F_n$ of $G$. Besides, since 
$\pi(G)$ has no invariant vector, there exists $g_0\in G$ such that 
$\pi(g_0)v\neq v$. Let $H$ be the subgroup of $G$ generated by the 
countable 
subset $\{g_0\}\cup\bigcup_n F_n$. Then the convex hull of $\alpha(H)(0)$ 
contains $D$. By Proposition \ref{cor:nilp_str_coh_triv}, $H$ has Property 
$H_{CT}$; it follows 
by Theorem \ref{thm:orbits_hct} that $D$ is parallel to the invariant 
vectors of $\pi(H)$, so that $v$ is contained in the 
$\pi(H)$-invariant vectors, 
a contradiction.\end{proof}

\begin{Cor}\label{densenilp}
Let $G$ be a topological group with Property $H_{CT}$. Let
$\mathcal{H}$ be a Hilbert space on which $G$ acts with
enveloping (respectively coarsely dense, resp. dense) image. Then
the action is by translations, defined by a continuous morphism:
$u:G\to (\mathcal{H},+)$ with enveloping (resp. coarsely dense,
resp. dense) image.\qed
\end{Cor}

\begin{Cor}
Let $G$ be a locally compact, compactly generated group with
Property $H_{CT}$, and let $\mathcal{H}$ be a (real) Hilbert
space. Then
\begin{itemize}
    \item $G$ has an isometric action on $\mathcal{H}$ with coarsely
dense (respectively enveloping) orbits if and only $\mathcal{H}$
has finite dimension $k$, and $G$ has a quotient isomorphic to
$\R^n\times\Z^m$, with $n+m\ge k$.
    \item $G$ has an isometric action on $\mathcal{H}$ with
dense orbits if and only $\mathcal{H}$ has finite dimension $k$,
and $G$ has a quotient isomorphic to $\R^n\times\Z^m$, with
$\max(n+m-1,n)\ge k$.
\end{itemize}
\label{cor:act_dense_nilp}
\end{Cor}

\begin{proof}
Let $\alpha$ be an affine isometric action of $G$ with enveloping
orbits (this encompasses all possible assumptions). By Corollary
\ref{densenilp}, the action is by translations; let $u$ be the
morphism $G\to(\mathcal{H},+)$; its image generates $\mathcal{H}$
as a topological vector space. Let $W$ denote the kernel of $u$.

Then $A=G/W$ is a locally compact, compactly generated abelian
group, which embeds continuously into a Hilbert space. By standard
structural results, $A$ has a compact subgroup $K$ such that $A/K$
is a Lie group. Since $K$ embeds into a Hilbert space, it is
necessarily trivial, so that $A$ is an abelian Lie
group without compact subgroup. 
Accordingly, $A$ is isomorphic to $\R^n\times\Z^m$ for some
integers $n,m$; the embedding of $A$ into $\mathcal{H}$ extends
canonically to a linear mapping of $\R^{n+m}$ into $\mathcal{H}$.
In particular $\mathcal{H}$ is finite-dimensional, of dimension
$k\le n+m$.

If the action has dense orbits, then either $m=0$ and $n\ge k$, or
$m\ge 1$ and $m\ge k-n+1$; this means that $k\le\max(n+m-1,n)$.
Conversely, if $k\le n+m-1$, then, since $\Z$ has a dense
embedding into the torus $\R^k/\Z^k$, $\Z^{k+1}$ has a dense
embedding into $\R^k$, and this embedding can be extended to
$\R^n\times\Z^m$.\end{proof}

From Proposition \ref{cor:nilp_str_coh_triv} and Corollary
\ref{cor:act_dense_nilp}, we deduce

\begin{Cor}\label{novirtuallynilpotent}
A compactly generated, nilpotent-by-compact locally compact group does not 
admit
any isometric action with enveloping (e.g. dense) orbits on an
infinite-dimensional Hilbert space.\qed
\end{Cor}


Proposition \ref{Prop:exist_metab_dense_Hilbert} on the one hand,
and Corollary \ref{novirtuallynilpotent} on the other, isolate the
first test-case for Navas' question:

\begin{quest}
    Does there exist a polycyclic group admitting
    an affine isometric action with dense
    orbits on an infinite-dimensional Hilbert space?
\end{quest}

Let us prove a related result for semisimple groups.

\begin{Thm}Let $G$ be a connected, semisimple Lie group. Then $G$
cannot act on a Hilbert space $\mathcal{H}\neq 0$ with coarsely
dense (e.g. dense) orbits.\label{thm:nodense_semisimple}
\end{Thm}
\begin{proof}Suppose by contradiction the existence of such an action
$\alpha$, and let $\pi$ denote its linear part. Then $\pi$ is
strongly cohomological. By Lemma \ref{trivialonZ(G)}, $\pi$ is
trivial on the centre of $G$. Thus the centre acts by
translations, generating a finite-dimensional subspace $V$ of
$\mathcal{H}$. The action induces a map $p:G\to V\rtimes\text{O}(V)$. 
Since $G$ is semisimple, the kernel of $p$ contains the sum
$G_{\text{nc}}$ of all noncompact factors of $G$, and thus factors
through the compact group $G/G_{\text{nc}}$. Thus $H^1(G,V)=0$, and
since $\pi$ is strongly cohomological, this implies that $V=0$.

It follows that $\alpha$ is trivial on the centre of $G$, so that
we can suppose that $G$ has trivial centre. Then $G$ is a direct
product of simple Lie groups with trivial centre. We can write
$G=H\times K$ where $K$ denotes the sum of all simple factors $S$
of $G$ such that $\alpha(S)(0)$ is bounded (in other words,
$H^1(S,\pi|_S)=0$). Then the restriction of $\alpha$ to $H$ also
has coarsely dense orbits. Moreover, every simple factor of $H$
acts in an unbounded way, so that, by a result of Shalom
\cite[Theorem 3.4]{ShaAff}\footnote{Shalom only states the result
for a simple group, but the proof generalizes immediately. See for
instance \cite{CoTeVa3} for another proof, based on the Howe-Moore
Property.}, the action of $H$ is proper. That is, the map $i:H\to
\mathcal{H}$ given by $i(h)=\alpha(h)(0)$ is metrically proper and
its image is coarsely dense. By metric properness, the subset
$X=i(H)\subset\mathcal{H}$ satisfies: $X$ is coarsely dense, and
every ball in $X$ (for the metric induced by $\mathcal{H}$) is
compact.

Suppose that $\mathcal{H}$ is infinite-dimensional and let us
deduce a contradiction. For some $d>0$, we have $d(x,X)\le d$ for
every $x\in\mathcal{H}$. If $\mathcal{H}$ is infinite-dimensional,
there exists, in a fixed ball of radius $7d$, infinitely many
pairwise disjoint balls $B(x_n,3d)$ of radius $3d$. Taking a point
in $X\cap B(x_n,2d)$ for every $n$, we obtain a closed, infinite
and bounded discrete subset of $X$, a contradiction.

Thus $\mathcal{H}$ is finite-dimensional; since every simple
factor of $H$ is non-compact, it has no non-trivial finite-dimensional 
orthogonal representation, so that the action is by
translations, and hence is trivial, so that finally
$\mathcal{H}=\{0\}$.\end{proof}

\begin{Rem}1) The same argument shows that a semisimple, linear algebraic
group over any local field, cannot act with coarsely dense orbits
on a Hilbert space.

2) The argument fails to work with enveloping orbits: indeed, in
$\ell^2_\R(\N)$, let $X$ denote the set sequences $(x_n)$ such
that $x_n\in 2^n\Z$ for every $n\in\N$. Then $X$ is coarsely dense
in $\ell^2_\R(\N)$, but, for the metric induced by $\mathcal{H}$,
every ball in $X$ is finite, hence compact. We do not know if a
semisimple Lie group (e.g. $\textnormal{SL}_2(\R)$) can act
isometrically on a nonzero Hilbert space with enveloping orbits.
\end{Rem}

\bigskip
\footnotesize
\noindent Yves de Cornulier\\
\'Ecole Polytechnique Fédérale de Lausanne (EPFL)\\
Institut de Géométrie, Algèbre et Topologie (IGAT)\\
1015 Lausanne, Switzerland\\
E-mail: \url{decornul@clipper.ens.fr}\\

\medskip

\noindent Romain Tessera\\
\'Equipe Analyse, Géométrie et Modélisation\\
Université de Cergy-Pontoise, Site de Saint-Martin\\
2 rue Adolphe Chauvin, 95302 Cergy-Pontoise Cedex, France\\
E-mail: \url{tessera@clipper.ens.fr}\\

\medskip

\noindent Alain Valette\\
Institut de Mathématiques - Université de Neuchâtel\\
Rue \'Emile Argand 11~- BP 158, 2009 Neuchâtel, Switzerland\\
E-mail: \url{alain.valette@unine.ch}


\begin{thebibliography}{KM98b}


\bibitem[CTV]{CoTeVa1} Yves {\sc de Cornulier}, Romain {\sc
Tessera}, Alain {\sc Valette}.
\newblock {\em Isometric group actions on Hilbert spaces: growth of
cocycles}. \newblock Preprint, 2005.

\bibitem[CLTV]{CoTeVa3} Yves {\sc de Cornulier}, Nicolas {\sc Louvet}, 
Romain {\sc
Tessera}, Alain {\sc Valette}.
\newblock {\em Howe-Moore Property and isometric actions on Hilbert
spaces}. \newblock In preparation, 2006.

\bibitem[Fe]{Fell} James M.G. {\sc Fell}.
\newblock {\em Weak containment and induced representations of groups}.
\newblock Canad. J. Math. {\bf 14}, 237--268, 1962.

\bibitem[Gu1]{Gui72} Alain {\sc Guichardet}.
\newblock {\em Sur la cohomologie des groupes topologiques II}.
\newblock Bull. Sci. Math. {\bf 96}, 305--332, 1972.

\bibitem[Gu2]{Guich} Alain {\sc Guichardet}. \newblock {\em Cohomologie
des groupes topologiques et des algèbres de Lie}. \newblock Paris,
Cédic-Nathan, 1980.



\bibitem[HiKa]{HiKa} Nigel {\sc Higson}, Gennadi {\sc Kasparov}.
\newblock {\em Operator K-theory for
groups which act properly and isometrically on Hilbert space}.
\newblock Electron. Res. Announc. Amer. Math. Soc. {\bf 3}, 
131--142, 1997.

\bibitem[MaVa]{MarVal} Florian {\sc Martin}, Alain {\sc Valette}.
\newblock {\em Reduced 1-cohomology of representation, the first
$\ell^p$-cohomology, and applications}. \newblock Preprint, 2005.

\bibitem[Sh1]{ShaAff} Yehuda {\sc Shalom}. \newblock {\em Rigidity, unitary representations
of semisimple groups, and fundamental groups of manifolds with
rank one transformation group}. \newblock Ann. of Math. (2) {\bf
152}(1), 113--182, 2000.

\bibitem[Sh2]{ShaInv} Yehuda {\sc Shalom}. \newblock {\em Rigidity of
commensurators and irreducible lattices}.
\newblock Invent. Math. {\bf 141}, 1--54, 2000.

\bibitem[Sh3]{Shalom04} Yehuda {\sc Shalom}. \newblock {\em Harmonic analysis,
cohomology, and the large scale geometry of amenable groups}.
\newblock Acta Math. {\bf 193}, 119--185, 2004.

\end{thebibliography}
\end{document}